\documentclass{ifacconf}

\usepackage{natbib}        

\usepackage{color,amsmath,amssymb,amsfonts,graphicx,url}  

\newcommand{\ml}{\left[ \begin{array}{cccccccccccccccccccccccccc}}
\newcommand{\mr}{\end{array} \right]}
\newcommand{\mll}{\left[ \begin{array}{c|c|c|c|c}}

\newcommand{\ii}{I\negmedspace I}
\newcommand{\xmin}{x_{\text{min}}}
\newcommand{\xmax}{x_{\text{max}}}

\newcommand{\rhomin}{\rho_{\min}}
\newcommand{\rhomax}{\rho_{\max}}
\newcommand{\rhominhat}{\hat{\rho}_{\min}}
\newcommand{\rhomaxhat}{\hat{\rho}_{\max}}

\newcommand{\supp}{\text{supp}}
\newcommand{\QED}{\begin{flushright} \qed \end{flushright}}

\newcommand{\subscr}[2]{#1_{\textup{#2}}}

\newcommand\oprocendsymbol{\hbox{$\square$}}
\newcommand\oprocend{\relax\ifmmode\else\unskip\hfill\fi\oprocendsymbol}

\renewcommand{\theenumi}{(\roman{enumi})}

\newcommand{\real}{\mathbb{R}}

\newtheorem{mydef}{Definition}[section]
\newtheorem{assumption}{Assumption}[section]

\begin{document}
\begin{frontmatter}

\title{Eulerian Opinion Dynamics with Bounded Confidence and Exogenous
  Inputs\thanksref{footnoteinfo}} 

\thanks[footnoteinfo]{This work was supported in part by the UCSB Institute for
    Collaborative Biotechnology through grant W911NF-09-D-0001 from the
    U.S.\ Army Research Office.}

\author[First]{Anahita~Mirtabatabaei \hspace{0.5in} Peng~Jia \hspace{0.5in} Francesco~Bullo} 

\address[First]{Center for Control, Dynamical Systems, and Computation, University of California, Santa Barbara, Santa Barbara, CA 93106, USA, (e-mail:  \{mirtabatabaei,pjia,bullo\}@engineering.ucsb.edu)}

\begin{abstract}                
The formation of opinions in a large population is governed by endogenous
(human interactions) and exogenous (media influence) factors.  In the
analysis of opinion evolution in a large population, decision making rules
can be approximated with non-Bayesian "rule of thumb" methods.  This paper
focuses on an Eulerian bounded-confidence model of opinion dynamics with a
potential time-varying input.  First, we prove some properties
of this system's dynamics with time-varying input. Second, we derive a
simple sufficient condition for opinion consensus, and prove the
convergence of the population's distribution with no input to a sum of
Dirac Delta functions. Finally, we define an input's \emph{attraction
  range}, and for a normally distributed input and uniformly distributed
initial population, we conjecture that the length of attraction range is an
increasing affine function of population's confidence bound and input's
variance.
\end{abstract}


\end{frontmatter}

\section{Introduction}

Decision-making in a society is a complex process determined by
\emph{endogenous} and \emph{exogenous} factors. The interaction of people
via in person meetings or online social networks is an endogenous factor.
One of the most influential exogenous factors is the mainstream media that
acts as a real-time input owing to its easy access to the public.
The early references~\cite{MS:61,SC-ES:77,CJ-EJ:86} propose models for
``continuous opinion dynamics,'' where opinions are represented by real
positive numbers. Such a real positive number represents the
worthiness of a choice or
the probability of choosing one decision over another.  A popular opinion
update rule is the non-Bayesian "rule of thumb" method of averaging
neighbors' opinions. This method provides a good approximation to the
behavior of a large population without relying on detailed social
psychological findings, see \cite{DA-AO:11}.
In our investigation, neighboring relation are defined based on
\emph{bounded confidence} (\cite{UK:00}), which means that an individual
only interacts with those whose opinions are close enough to its own.  This
idea reflects: 1) \emph{filter bubbles}, a phenomenon in which websites use
algorithms to show users only information that agrees with their past
viewpoints (\cite{EP:11}); and 2) \emph{selective exposure}, a
psychological concept broadly defined as "behaviors that bring the
communication content within reach of one's sensory apparatus"
(\cite{DZ-JB:85,WM-FC-RS:08}).
According to \cite{CC-FF-PT:08,JMH:08}, models of opinion dynamics can be
described by either a \emph{Lagrangian} or an \emph{Eulerian} method.  A
Lagrangian description focuses on changes in each agent's opinion; an
Eulerian description focuses on the changes in agents population in one
opinion interval as time progresses.  A Lagrangian model of opinion
dynamics is defined over a continuous (\cite{VDB-JMH-JNT:09}) or discrete
state space (\cite{RH-UK:02,VDB-JMH-JNT:09,AO-JNT:07,JL:09}) if the number
of agents are infinite or finite, respectively.  An Eulerian model of
opinion dynamics is defined over a continuous
(\cite{JL:08,CC-FF-PT:08,GC-FF:09}) or discrete state space
(\cite{SB-RL-TA:11}) depending on whether the opinion grid size is
converging to zero or not, respectively.

\cite{RH-UK:02} formulated a Lagrangian bounded confidence model, referred
to as the HK model, where agents synchronously update their opinions.
Accordingly, an Eulerian HK model is defined over a continuous state space
(\cite{JL:07b,CC-FF-PT:08}), where a \emph{mass distribution} over an
opinion set is being updated by a \emph{flow map} under the influence of an
exogenous input.  Previously, \cite{CC-FF-PT:08} proved the convergence of
a variation of Eulerian HK model both in discrete and continuous time. In
their model, the weights that two opinion values assign to each other are
equal, and this symmetry preserves the global average during the
evolution. In this context, we consider a more general Eulerian HK model
where the symmetric weight constraint is relaxed. Specifically, the weight
an opinion assigns to other opinions is a function of the integral of the
mass distribution (and of the exogenous input measure) in that opinion's
confidence bound. Since the measures on different opinions' confidence
bounds are not necessarily equal, the weights assigned to different
opinions are generally asymmetric, and thus the global average is not
preserved.

Media influences decisions by employing some well known techniques such as
repeated exposure to experts' messages (\cite{YL-JB-DL:11}). Increasingly, however, the
message sent by the media is restated by public blogs with bias.
Reflecting these properties, we model the media exogenous input in our
model as a background Gaussian signal centered at the opinion of an expert.
Since the media influence on the public depends on public's attitude
towards it (\cite{JML:07}), we assume that each agent associated with one
opinion receives exogenous input information only within its opinion
confidence ranges.  The effect of media on opinion formation with pairwise
gossip interactions is numerically analyzed in \cite{LB-RM-FS:10}.

The contributions of this paper can be summarized as follows.  1) We
propose a reasonable model for exogenous inputs in the Eulerian HK
opinion-dynamics model.   We derive a simple sufficient condition for the
system to reach opinion consensus.
2) We establish some important properties of the Eulerian HK model with a
time-varying input.  Under mild technical assumptions (the initial opinion
is a finite and absolutely continuous
mass distribution over the opinion set), we show that the opinion update
via an Eulerian flow map has the following properties: i) the mass
distribution on opinions remains finite and absolutely continuous; ii) the
flow map preserves opinion order, due to the homogeneity of confidence
bounds; and iii) the flow map is bi-Lipschitz.
This analysis also leads to a convergence proof of the mass distribution to
a sum of Dirac Delta functions.
3) We represent the exogenous input by a background Gaussian distribution centered at the advertised opinion.
 We introduce the \emph{attraction range} of an input, which is the largest range of opinions that the input can attract to its center. We conjecture a linear relation between attraction range, input's variance, and confidence bound. 
 Accordingly, we compare two different manipulation strategies that aim to increase the population who vote positively in finite time.

\section{Mathematical Model}\label{SecModel}

%
In this paper we describe the process of opinion exchange in a large population via a sequence of finite Borel measures. 
This approach is inspired by \cite{CC-FF-PT:08}, where the \emph{mass distribution} of agents over the opinion set  is represented by $\mu_t: \real \rightarrow \real_{\ge 0}$ at discrete time steps $t$, 
whose sum over opinion space is preserved over time.
 The opinion set belongs to a continuous state space in $\real$, and each opinion value is denoted by independent variable $x$.
In view of state space's continuity, the mass distribution $\mu_t(x)$ is assumed to be a finite Borel probability measure on $R$. 
The value $\mu_t(dz)$ at $x$, denoted by $d\mu_t(x)$, represents the infinitesimal population whose opinion is equal to $x$ at time $t$.
At $t+1$, this population updates its opinion to $\gamma_t(x)$, defined as the \emph{flow map} of mass distribution  $\gamma_t: \supp \mu_t \subseteq \real \rightarrow \real$, where $\supp \mu_t$  denotes the \emph{support} of the measure $\mu_t$, that is, the set of all points $x \in R$ for which every open neighborhood of $x$ has positive measure.
Here, the flow map is defined in compliance with the Lagrangian HK rules,
\begin{equation}\label{flowmap}
\gamma_t(x) = \frac{\int_{[x-r,x+r]} z d\mu_t(z) + \int_{[x-r,x+r]} z du_t(z) }{\int_{[x-r,x+r]} d\mu_t(z) + \int_{[x-r,x+r]} du_t(z) },
\end{equation}
where $r$ is the confidence bound of agents, and $u_t$ is the distribution of exogenous background input at time $t$, which is also assumed to be a finite Borel probability measure for simplicity of analysis. 
Now, the mass distribution can be tracked by the recurrence relation $\mu_{t+1} = \gamma_t \# \mu_t,$
where $\gamma_t \# $ denotes the \emph{push forward} of a measure via the flow map $\gamma_t$ (\cite{CC-FF-PT:08}). Moreover, for every Borel set $E \in \real$,
\begin{equation}\label{eqpushfw}
\mu_{t+1}(E) = \mu_t(\gamma_t^{-1}(E)),
\end{equation}
where $\gamma_t^{-1}(E)$ is the preimage of set $E$ under flow map $\gamma_t$, that is not necessarily invertible. 
\begin{mydef}[Eulerian HK System with Input]
We call the dynamical system in which a mass distribution $\mu_t$ defined over a continuous state space is being pushed forward with flow map~\eqref{flowmap} under the influence of input $u_t$ a \emph{Eulerian HK system with Input}.
\end{mydef}

In the remainder of this section, we introduce a lemma that defines a bound on the value of flow map, then we present main assumptions of this paper.
From here on, absolute continuity of any measure $\mu$ with respect to Lebesgue measure $\mathcal{L}^1$ is denoted with $\mu \ll \mathcal{L}^1$; the smallest and largest opinions along $\supp \mu_t$ are denoted by $\xmin(t)$ and $\xmax(t)$, respectively; and the length of interval $\xmax(t) - \xmin(t)$ is denoted by $|\supp \mu_t|$.  Moreover, flow map  $\gamma_t(x)$ is called bi-Lipschitz, if for any $x,y \in \supp \mu_t$ there exists $L_t \ge 1$ such that 
\begin{equation}\label{EqBilip}
|y - x|/L_t \le |\gamma_t(y) - \gamma_t(x) | \le L_t |y-x|.
\end{equation}

\begin{lem}[Bounds on mass average]\label{Lemboundedave}
For a finite mass distribution $\mu \ll \mathcal{L}^1$, assume that its support is a closed bounded interval of $R$, and its density function $\rho_{\mu}(x)\geq 0$, satisfies $\rho_{\mu}(x) \in [\rhomin,\rhomax]$ for all $x \in \supp \mu$, where $0 < \rhomin\le\rhomax<\infty$. Then, for all $a,b \in \supp \mu$, the opinion average over $[a,b]$, denoted by $A$,
$$ A:= \frac{\int_{[a,b]} z d\mu(z)}{\int_{[a,b]} d\mu(z)}=\frac{\int_a^b z \rho_{\mu}(z) d z}{\int_a^b \rho_{\mu}(z) d z}$$
 can be bounded as follows
\begin{equation}\label{strictineq}
a \;  < \; \frac{ b + a \sqrt{\rhomax/\rhomin} }{1 +  \sqrt{\rhomax/\rhomin}} \; 
 \le
\; A
  \; \le \; \frac{ a + b \sqrt{\rhomax/\rhomin} }{1 +  \sqrt{\rhomax/\rhomin}} \; < \; b.
\end{equation}
\end{lem}

\begin{pf} Since $\mu$ is assumed to be a finite absolutely continuous measure, there exists a Lebesgue integrable density function $\rho_{\mu}$ such that $\mu(E)=\int_{E}\rho_{\mu}(z)dz$ for all Borel subsets $E \in \real$. 
Here, we prove the upper bound of the average, and proof to the lower bound is similar.
We maximize $A$ for the following \emph{step density function} over the variable $c \in [a, b]$:
\begin{equation*}
\rho_s(x)= 
\begin{cases}
 \rhomin, & \mbox{if }  x \in [a, c),\\
 \rhomax,  & \mbox{if }  x \in [c, b],
\end{cases} 
\end{equation*}
According to the first mean value theorem for integrals, one can show that for any bounded density function $\rho: [a,b] \rightarrow [\rhomin,\rhomax]$, there exists $c \in [a, b]$ such that the averages of $\rho$ and $\rho_s$ over $[a,b]$ are equal. For $\rho_s$ we have
\begin{equation*}
A = \frac{(c^2-a^2)\rhomin/2 +(b^2-c^2)\rhomax/2 }{(c-a)\rhomin +(b-c)\rhomax} =: \frac{f}{g}.
\end{equation*}
Owing to the differentiability of $A$ with respect to $c$, the maximum of $A$ over $c$ can be computed by letting $\partial A /\partial c$ equal to zero.
\begin{equation*}
\frac{\partial A}{\partial c} = \frac{c (\rhomin - \rhomax)g -  (\rhomin - \rhomax)f }{g^2} = 0.
\end{equation*}
Thus, the critical point $c = f/g = A$ gives maximum $A$,
\begin{multline*}
 2(c^2 -ac)\rhomin + 2(bc - c^2)\rhomax 
 \\ =(c^2 - a^2)\rhomin + (b^2 - c^2)\rhomax 
\Rightarrow \; c = \frac{ a + b \sqrt{\rhomax/\rhomin} }{1 +  \sqrt{\rhomax/\rhomin}}.
\end{multline*} \QED
\end{pf}
\begin{assumption}\label{assumptions}
For an Eulerian HK system with input, $\mu_0 \ll \mathcal{L}^1$ is finite and $\supp \mu_0$ is a closed bounded interval, and $u_t \ll \mathcal{L}^1$ for all $t \ge 0$.
\end{assumption}
\begin{assumption}\label{assumSupU}
The set $\supp u_t$ is contained in the set $\supp \mu_t$ for all $t \ge 0$. 
\end{assumption}
The interpretation of Assumption~\ref{assumSupU} is that the manipulator can only advertise for opinions with non-zero population assigned to them. In other words, the manipulator disregards the opinions that no body believes in.

 \section{Dynamics of the Model}\label{SecProp}

This section first analyzes some fundamental properties of Eulerian HK systems with time-varying inputs. Employing these properties, the convergence of Eulerian HK systems with no input is established. The following lemma is a simple sufficient condition for consensus whose proof is omitted in the interest of brevity.

\begin{lem}[Sufficient condition for consensus]\label{lemmaConsensus}
In an Eulerian HK system with input $u_t$, assume that $u_t$ satisfies Assumption~\ref{assumSupU} and $\mu_0$ is a finite measure with closed bounded support. If  $\mu_0$ and $u_t$ are distributed symmetrically around the center of $\supp \mu_0$, and $|\supp \mu_0| < 2r$, then the mass distribution reaches an opinion consensus in finite time. 
\end{lem}
\begin{thm}[Properties of Eulerian HK system with input]\label{propertiesProp}
If an Eulerian HK system with input satisfies Assumption~\ref{assumptions}, then for all $t \ge 0$ such that $|\supp \mu_t| > 2r$ and $\supp u_t \cup \supp \mu_t$ is a closed bounded interval,
\begin{enumerate}
\item\label{Prop1bndMu} $\mu_t \ll \mathcal{L}^1$ is finite and $\supp \mu_t$ is a closed interval;
\item\label{Prop1monotone} for any $x,y \in \supp \mu_t$, if $x < y$, then $\gamma_t(x) < \gamma_t(y)$; and
\item\label{Prop1lip} $\gamma_t(x)$  is bi-Lipschitz with respect to $x$.
\end{enumerate}
\end{thm}
\begin{pf} Here, we first prove that if statement~\ref{Prop1bndMu} holds at any time $t$, then statements~\ref{Prop1monotone} and \ref{Prop1lip} will hold at $t$. Next, if the three statements hold at any $t$, then statement~\ref{Prop1bndMu} holds at $t+1$. Finally, since $\mu_0$ satisfies statement~\ref{Prop1bndMu}, the three statements hold for all $t$. For brevity, we denote the sum of the mass and input distributions with $\nu_t := \mu_t + u_t$.
Since $u_t$ satisfies Assumption~\ref{assumptions}, if statement~\ref{Prop1bndMu} holds at any $t$, then 
$\nu_t \ll \mathcal{L}^1$ is finite and $\supp \nu_t$ is a closed bounded interval. Hence, 
$\nu_t$'s density function $\rho_t(x)\geq 0$ exists and satisfies $\rho_t(x) \in [\rhomin(t),\rhomax(t)]$ for all $x \in \supp \nu_t$ with $0 < \rhomin(t) \le \rhomax(t) <\infty$.

Regarding part~\ref{Prop1monotone}, for any $x, y \in \supp \mu_t$ and $x< y$, since $x \pm r$ or $y \pm r$ may not belong to $\supp \nu_t$, 
\begin{equation}\label{eqgammaxy}
\gamma_t(x) = \frac{\int_a^b z \rho_t(z) dz}{\int_a^b \rho_t(z) dz}, \quad \gamma_t(y) = \frac{\int_p^q z \rho_t(z) dz}{\int_p^q \rho_t(z) dz},
\end{equation}
where $[a,b] = [x-r,x+r] \cap \supp \nu_t$ and $[p,q] = [y-r,y+r] \cap \supp \nu_t$. Equivalently,
\begin{equation}\label{eqflowx}
\gamma_t(x) = \frac{\int_{a}^{p} z \rho_t(z) dz + \int_{p}^{b} z \rho_t(z) dz }{\int_{a}^{p} \rho_t(z) dz + \int_{p}^{b} \rho_t(z) dz} =: \frac{\hat{S}_1 + \hat{S}_2}{S_1+S_2},
\end{equation}
\begin{equation}\label{eqflowy}
\gamma_t(y) = \frac{\int_{p}^{b} z \rho_t(z) dz + \int_{b}^{q} z \rho_t(z) dz }{\int_{p}^{b} \rho_t(z) dz + \int_{b}^{q} \rho_t(z) dz} =: \frac{\hat{S}_2 + \hat{S}_3}{S_2+S_3}. 
\end{equation}
It follows from properties of $\nu_t$ that Lemma~\ref{Lemboundedave} holds, and considering the integration intervals of $\hat{S}_i$'s and $S_i$'s, for nonzero $S_i$'s and any $i < j$, we have
$\hat{S}_i/S_i < \hat{S}_j/S_j$, and thus $\hat{S}_iS_j<\hat{S}_jS_i.$
Notice that based on assumption $|\supp \mu_t| > 2r$, at least one of the $S_1$ or $S_3$ should be nonzero, moreover, since $\supp \nu_t$ is a closed interval, the terms $S_1+S_2$ and $S_1+S_3$ are nonzero. Consequently, only one term out of the three terms $S_1$, $S_2$ and $S_3$ can be zero, and the following inequality always holds:
\begin{multline*}
\hat{S}_1S_2 + \hat{S}_1S_3 + \hat{S}_2S_2 +\hat{S}_2S_3 < \hat{S}_2S_1 + \hat{S}_2S_2 + \hat{S}_3S_1 +\hat{S}_3S_2, \\
\Rightarrow \;\;  \frac{\hat{S}_1 + \hat{S}_2}{S_1+S_2} 
< \frac{\hat{S}_2 + \hat{S}_3}{S_2+S_3} \;\;
 \Rightarrow \;\;  \gamma_t(x) < \gamma_t(y).
\end{multline*}

Regarding part~\ref{Prop1lip}, the bi-Lipschitz property of $\gamma_t(x)$
asserts that for any $x,y \in \supp \mu_t$ equation~\eqref{EqBilip} holds for some $L_t \ge 1$.
Assume that $x < y$, and by part~\ref{Prop1monotone}, $\gamma_t(x) < \gamma_t(y)$. Then, two different cases are possible:

1) $y - x \ge 2r$, hence,
$ \gamma_t(y) - \gamma_t(x) < y - x + 2r \le 2(y -x),$
and it follows from Lemma~\ref{Lemboundedave} that
$$ \gamma_t(y) - \gamma_t(x) > \frac{(b-a) + (q - p)}{1+\sqrt{\rhomax(t)/\rhomin(t)}},$$
where the flow maps are given by equations~\eqref{eqgammaxy}.
Since $y -x \le |\supp \mu_t|$,
\begin{multline*} 
\frac{(b-a) + (q - p)}{1+\sqrt{\rhomax(t)/\rhomin(t)}}  
=  \frac{(b-a + q - p) (y -x)}{(1+\sqrt{\rhomax(t)/\rhomin(t)}) (y -x)} \\
 \ge \frac{(b-a + q - p) (y -x)}{|\supp \mu_t|(1+\sqrt{\rhomax(t)/\rhomin(t)})}.
\end{multline*}
Finally, $$L_t = \min \{2, \frac{|\supp \mu_t|(1+\sqrt{\rhomax(t)/\rhomin(t)})}{b-a + q - p} \}.$$

2) $y - x < 2r$, by equations~\eqref{eqflowx} and \eqref{eqflowy}, $\gamma_t(y) - \gamma_t(x) $ can be written as a function of $S_i$'s and $\hat{S}_i$'s and be upper-bounded as follows
\begin{multline*} 
 \gamma_t(y) - \gamma_t(x) 
 < \big(2r |q|\rhomax(t)^2(y-x) +|q|\rhomax(t)^2(y-x)^2 \\+ 2r |q|\rhomax(t)^2(y-x)\big)/r^2\rhomin(t)^2.
\end{multline*}
Again since $y-x < 2r$ and $|q| \le \max\{|\xmax(t)|, |\xmin(t)|\}$,
\begin{equation*}
\gamma_t(y) - \gamma_t(x)  <  \frac{ 6r \max\{|\xmax(t)|, |\xmin(t)|\} \rhomax(t)^2}{r^2\rhomin(t)^2} (y-x),
\end{equation*}
where we denote the coefficient of $y -x$ by $L_1$.
As stated above, either $S_1$ or $S_3$ is nonzero. Without loss of generality, assume that $S_1$ is nonzero, and hence $p-a = y -x$. It follows from 
$(\hat{S}_3 + \hat{S}_2)(S_3+S_2) \ge \hat{S}_2/S_2$ that
\begin{equation}\label{gammaineq}
\gamma_t(y) - \gamma_t(x)  \ge  \frac{\hat{S}_2}{S_2} - \frac{\hat{S}_1 + \hat{S}_2}{S_1+S_2} =: c(x_2-x_1),
\end{equation}
where $x_1 = \hat{S}_1/S_1, x_2 = \hat{S}_2/S_2,$ and
\begin{multline*}
c = \frac{S_1}{S_1+S_2} \ge \frac{(p-a)\rhomin(t)}{(b-p)\rhomax(t)+(p-a)\rhomax(t)},\\
\Rightarrow \gamma_t(y) - \gamma_t(x)  \ge 
\frac{(p-a)\rhomin(t)}{(b-a)\rhomax(t)} (x_2 -x_1) \\
> \frac{(y-x)\rhomin(t)}{2r\rhomax(t)} (x_2 -x_1).
\end{multline*}
By Lemma~\ref{Lemboundedave}, 
\begin{equation*}
x_2 - x_1 > \frac{b-p + p-a}{\sqrt{\rhomax(t)/\rhomin(t)} +1} 
\ge \frac{r}{\sqrt{\rhomax(t)/\rhomin(t)} +1} . 
\end{equation*}
Therefore, 
\begin{equation*}
\gamma_t(y) - \gamma_t(x)  > 
\frac{\rhomin(t)}{2\rhomax(t)(\sqrt{\rhomax(t)/\rhomin(t)} +1)} (y-x), 
\end{equation*}
denoting coefficient of $y -x$ by $L_2$, $L_t = \min \{L_1, 1/L_2\}$.

Regarding part~\ref{Prop1bndMu}, we now prove that if statements~\ref{Prop1bndMu}, \ref{Prop1monotone}, and \ref{Prop1lip} hold at time $t$, then statement~\ref{Prop1bndMu} holds at $t+1$.
First, we prove that the flow map $$\gamma_t(x) = \frac{\int_{x-r}^{x+r} z \rho_t(z) dz}{\int_{x-r}^{x+r} \rho_t(z) dz} =: \frac{f(x)}{g(x)}$$ 
is continuous.
Knowing that if two functions $f$ and $g$ are continuous and $g\neq 0$, then $f/g$ is also continuous, we show the continuity of the function $f(x)$ at all points $c \in \supp \mu_t$, and the proof to the continuity of the denominator is similar. For all $x \in \supp \mu_t$, $g(x) > 0$, and
\begin{multline*}
\lim_{x \rightarrow c} f(x) = \lim_{x \rightarrow c} \int_{x-r}^{x+r} z d\rho(z)= \lim_{\pm\epsilon \rightarrow 0} \int_{c\pm\epsilon-r}^{c\pm\epsilon+r} z d\rho(z)\\
 =  \lim_{\pm\epsilon \rightarrow 0} \big( -\int_{c-r}^{c\pm\epsilon-r} z d\rho(z) + \int_{c+r}^{c\pm\epsilon+r} z d\rho(z) \big) + \int_{c-r}^{c+r} z d\rho(z) \\
 = \lim_{\pm\epsilon \rightarrow 0} \big( \pm\epsilon (r - c) \rho(c-r) +\pm\epsilon (c+r) \rho(c+r) \big) + f(c).
\end{multline*}
Due to the finiteness and absolute continuity of $\nu_t$, $f(c)$ exists for all $c\in \supp \mu_t$ and above limit converges to zero, hence, $\lim_{x \rightarrow c} f(x) = f(c)$. 
We have shown that $\gamma_t(x)$ is strictly monotone and continuous with respect to $x$, therefore, this map is also invertible.
Second, we prove absolute continuity of $\mu_{t+1}$. It is shown that $\gamma_t$ has the following properties:
1) since any continuous function defined on Borel sets is a Borel measurable function, $\gamma_t^{-1}(E)$ is Borel measurable for any Borel set $E \in \real$. 
2) the bi-Lipschitz map $\gamma_t$ satisfies 
$\mathcal{L}(\gamma_t^{-1}(E)) \le C_t \mathcal{L}(E)$ for some constant $C_t \in \real_{>0}$. According to Theorem 2 in \cite{BP-AT:10}, if the flow map $\gamma_t$ satisfies above two properties and $\mu_t \ll \mathcal{L}^1$, then  $\mu_{t+1}  \ll \mathcal{L}^1$.
Third, a continuous function maps a compact set to another compact set, hence, $\gamma_t$ maps the closed bounded interval $\supp \mu_t$ to another closed bounded interval $\supp \mu_{t+1}$.
Fourth, we establish bounds on $\mu_{t+1}$'s density function. For any $x , y \in \supp \mu_{t+1}$ and $x<y$, equation~\eqref{eqpushfw} gives
\begin{equation*}
\int_x^y \hat{\rho}_{t+1}(z)dz = \int_{\gamma_t^{-1}(x)}^{\gamma_t^{-1}(y)} \hat{\rho}_t(z) dz,
\end{equation*}
where $\hat{\rho}_{\tau}(z)$ is $\mu_{\tau}$'s density function, and $ 0 < \rhominhat(\tau) \le \hat{\rho}_\tau(z)\le \rhomaxhat(\tau) < \infty$ over $\supp \mu_\tau$. Therefore,
\begin{multline*}
 (\gamma_t^{-1}(y) - \gamma_t^{-1}(x)) \rhominhat(t) \le \int_{\gamma_t^{-1}(x)}^{\gamma_t^{-1}(y)} \hat{\rho}_t(z) dz 
 \\\le  (\gamma_t^{-1}(y) - \gamma_t^{-1}(x)) \rhomaxhat(t),\\
\Rightarrow\; \frac{y-x}{L_t} \rhominhat(t) \le \int_{\gamma_t^{-1}(x)}^{\gamma_t^{-1}(y)} \hat{\rho}_t(z) dz \le  L_t(y-x) \rhomaxhat(t), \\
\Rightarrow\; \frac{y-x}{L_t} \rhominhat(t) \le \int_x^y  \hat{\rho}_{t+1}(z)dz  \le  L_t(y-x) \rhomaxhat(t).
\end{multline*}
The limit of above inequality as $y$ converges to $x$ gives
\begin{equation*}
\frac{1}{L_t} \rhominhat(t) \le \hat{\rho}_{t+1}(x) \le  L_t \rhomaxhat(t), \;\;\; \forall x \in \supp \mu_{t+1}.
\end{equation*}
Finally, since $\supp \mu_t$ is bounded for all $t\ge 0$, we have $\mu_{t+1}(\supp \mu_{t+1}) = \mu_t(\gamma_t^{-1}(\supp \mu_{t+1})) = \mu_t(\supp \mu_t)$. Therefore, if $\mu_t$ is finite, then $\mu_{t+1}$ is finite. \hspace{.5in} \qed
\end{pf}
The following lemmas are employed in the proof of Theorem~\ref{theoremConvergence}, and their proofs are omitted for brevity.
\begin{lem}\label{lemmaSuppShrink}
If an Eulerian HK system with input satisfies Assumptions~\ref{assumptions} and \ref{assumSupU}, then for all $t \ge 0$ such that $|\supp \mu_t| > 2r$, (i) $\supp \mu_t$ strictly contains $\supp \mu_{t+1}$, and (ii) $\xmin(t+1) = \gamma_t(\xmin(t))$,  and  $\xmax(t+1) = \gamma_t(\xmax(t))$.
\end{lem}
Notice that if $\supp u_t$ is not contained in $\supp \mu_t$, then $\supp \mu_{t+1}$ does not necessarily contain in $\supp \mu_t$.
\begin{lem}\label{LemmaMaximize}
Consider $x < y \in \real$ and $A, B \in \real_{>0}$, then 
$$ \max_{A \in [a_1,a_2], B \in [b_1,b_2]}\frac{Ax + By}{A+B} = \frac{a_1x+ b_2y}{a_1+b_2}.$$
\end{lem}
In the following context, 
we call a single point $x \in \real$ an atom w.r.t. a measure $\mu$ if $x\in \supp\mu$ and $\mu(x)>0$. 
Moreover, if every  $\mu$-measurable set of positive measure contains an atom, then $\mu$ is purely atomic or atomic in short.
\begin{thm}\label{theoremConvergence} 
Consider an Eulerian HK system with no input where the initial condition satisfies that  $\mu_0 \ll \mathcal{L}^1$ is finite and $\supp \mu_0$ is a closed bounded interval. If $|\supp \mu_t| > 2r$ for all $t\ge 0$, then $\mu_t$ converges in the weak-star topology to 
an atomic measure, whose atoms are a distance greater than $r$ apart from one another. 
\end{thm}
{\bf Sketch of proof.} This system satisfies conditions of Theorem~\ref{propertiesProp} and Lemma~\ref{lemmaSuppShrink}.
Therefore, $\mu_t  \ll \mathcal{L}^1$, $\supp \mu_t$ is a closed bounded interval, and $|\supp \mu_t|$ strictly decreases at each iteration. 
Since $\xmin(t)$ is strictly increasing  and
$\supp \mu_t$ is a subset of $\supp \mu_0$, there exists an opinion $x_1$ in the interior of $\supp \mu_t$ such that $\lim_{t\rightarrow \infty} \xmin(t) = x_1$.  Thus, there exists $\tau$ after which $x_1 - \xmin(t) < r$.
 
 First, we prove that the mass distribution over interval $(x_1,x_1+r)$ converges to zero.
 Let us denote  $\mu([\xmin(t),x_1])$ and $\mu((x_1,\xmin(t)+r])$ by $\hat{S}_t$ and $S_t$, respectively. Since $\xmin(t)$ converges to $x_1$, convergence of $S_t$ to zero and $\mu_t((x_1,x_1+r))$ to zero are equivalent. 
 Suppose by contradiction that there exists a constant $S_{\min} > 0$ such that $ S_t\ge S_{\min}$ for all $t \ge \tau$.
Denoting the density function of $\mu_t$ by $\rho_t$, we define
$$y_t := \frac{\int_{x_1}^{\xmin(t)+r} z \rho_t(z) dz}{\int_{x_1}^{\xmin(t)+r} \rho_t(z) dz},$$
 then Lemma~\ref{Lemboundedave} tells us that $y_t > x_1$ for all $t \ge \tau$.
Since $\xmin(t)$ converges to $x_1$ as time goes to infinity,
$$\lim_{t\rightarrow \infty} \int_{\xmin(t)}^{x_1} z \rho_t(z) dz = x_1\lim_{t\rightarrow \infty} \hat{S}_t.$$  
By Lemma~\ref{lemmaSuppShrink}, $\xmin(t+1) = \gamma_t(\xmin(t))$, hence
\begin{equation*}
 \lim_{t\rightarrow \infty}  \xmin(t+1) = \frac{ x_1\lim_{t\rightarrow \infty} \hat{S}_t+ y_tS_t}{\lim_{t\rightarrow \infty} \hat{S}_t+S_t},
\end{equation*}
and by Lemma~\ref{LemmaMaximize}
 \begin{equation*}
 \lim_{t\rightarrow \infty} \xmin(t+1) > \frac{x_1\lim_{t\rightarrow \infty} \hat{S}_t+ x_1S_{\min}}{\lim_{t\rightarrow \infty} \hat{S}_t+S_{\min}} = x_1
 \end{equation*}
which contradicts monotonic convergence of $\xmin(t)$ to $x_1$.  

Second, we show that $\hat{S}_t$ converges to an atomic measure. 
According to the first part of this proof, there exists a time $T \ge \tau$ after which $\mu_t((x_1,x_1+r)) < \epsilon$ for any $\epsilon \in \real_{>0}$.
Consider $x \in (\xmin(t), x_1)$ for which $\mu_t(\mathcal{B}_{\delta}(x)) > \epsilon$ for  $t \ge T$, where $\mathcal{B}_{\delta}(x)$ is an open ball centered at $x$ with infinitesimal radius $\delta \in \real_{>0}$.
Suppose by contradiction that $\gamma_t(x) > x_1$. Hence, $\mu_{t+1}(\mathcal{B}_{\delta}(\gamma_t(x))) > \epsilon$, while
according to Lemma~\ref{Lemboundedave}, $\gamma_t(x) < x_1 + r$.  In other words, there exists $y \in (x_1,x_1+1)$ such that $\mu_{t+1}(\mathcal{B}_{\delta}(y)) > \epsilon$, which contradicts the first part of this proof.
Therefore, for all $t \ge T$,  if $x \in (\xmin(t), x_1)$ and $\rho_t(x) > 0$, then $\gamma_t(x) \le x_1$, and thus $\hat{S}_t$ is a non-decreasing function of time.
Theorem~\ref{propertiesProp} tells us that $\mu_t$ is finite for all $t$, hence $\hat{S}_t$ is upper-bounded. Consequently, there exists $\hat{S}_1 > 0$ such that $ \lim_{t\rightarrow \infty} \hat{S}_t = \hat{S}_1$.

Third, owing to the convergence of $\mu_t((x_1,x_1+r))$ to zero, 
\begin{equation*}
\lim_{t \rightarrow \infty}\gamma_t(x_1+r) 
=  \lim_{t \rightarrow \infty}  \frac{\int_{x_1+r}^{x_1+2r} z d\mu_t(z) }{\int_{x_1+r}^{x_1+2r} d\mu_t(z) }  > x_1+r.
\end{equation*}
Consequently, the first and second part of this proof can be repeated for $x_1+r$ as the $\xmin(t)$ of the mass distribution in the interval $(x_1+r,\xmax(t)]$ for $t\ge T$ and so on.

Finally, for every bounded and continuous test function $\eta$
\begin{multline*}
\lim_{t \rightarrow \infty} \int_{\real} \eta(z)d \mu_t(z) = x_1\eta(x_1)\hat{S}_1  + x_2 \eta(x_2)\hat{S}_2 + \dots. 
\end{multline*}
 Hence, $\mu_t$ converges in weak-star topology to an atomic measure whose atoms are at least $r$ distance far apart.   \qed

\section{Discussion on Exogenous Event}\label{SecInput}

\begin{figure}\label{figLinearrelation}
  \includegraphics[width=.24\textwidth,keepaspectratio]{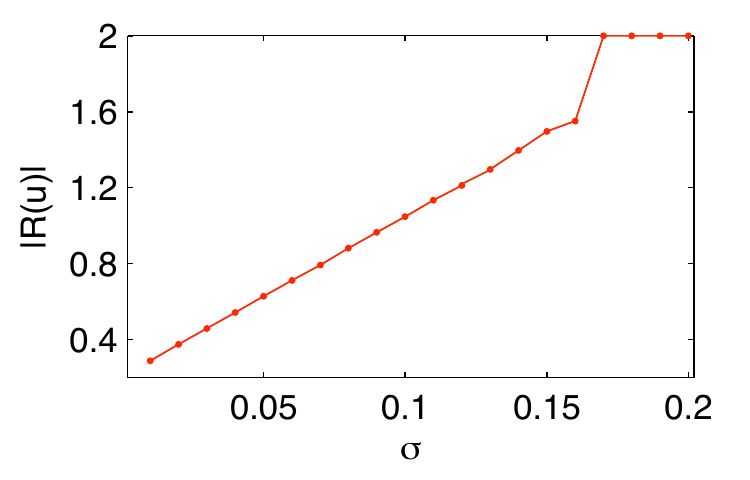}
  \includegraphics[width=.24\textwidth,keepaspectratio]{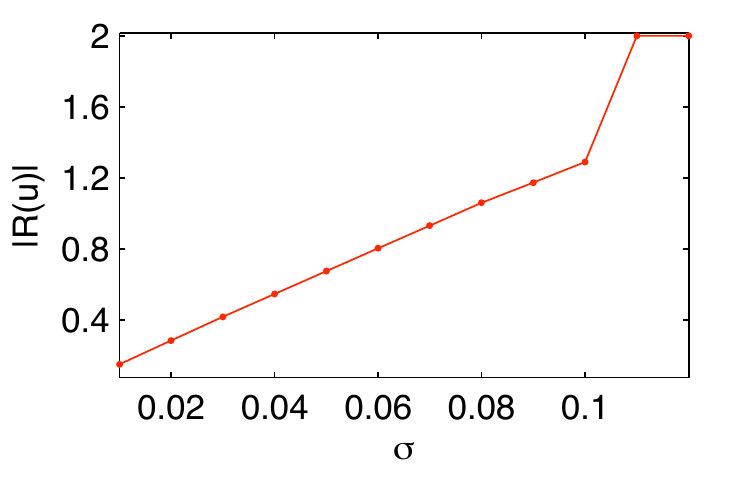}\\
 \includegraphics[width=.24\textwidth,keepaspectratio]{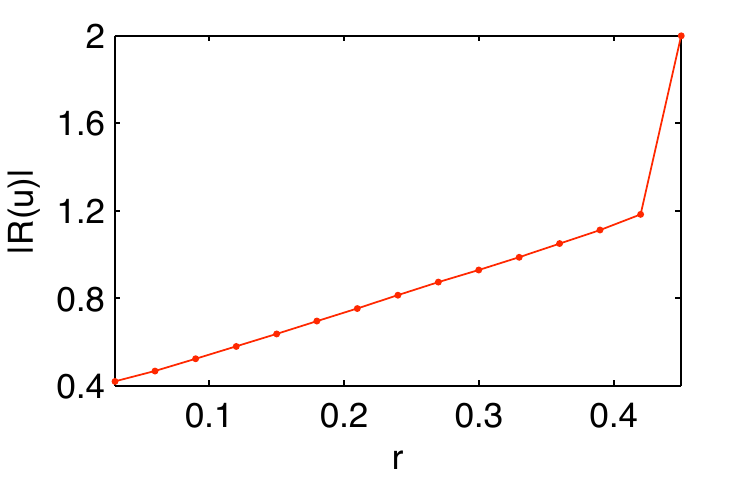}
  \includegraphics[width=.24\textwidth,keepaspectratio]{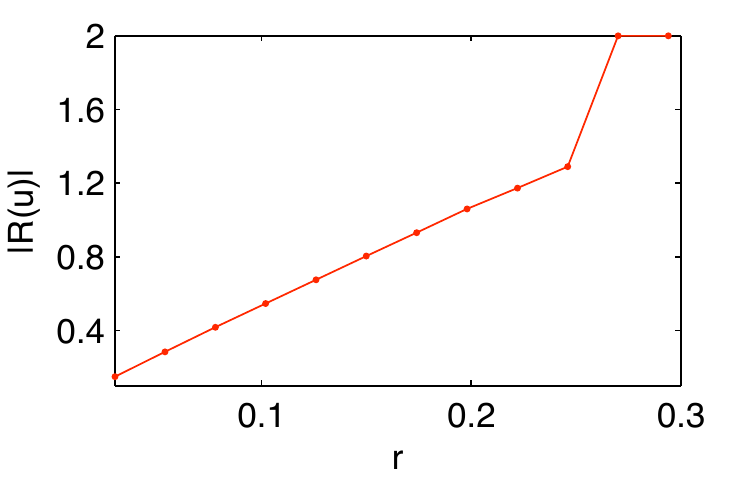}
\caption{In evolutions of Eulerian HK systems with uniform initial
  distribution $\mu_0 \sim \mathcal{U}(-x_0,x_0)$ and input $u \sim
  \subscr{\mathcal{N}}{truncated}(0,\sigma^2)$, $|\mathcal{R}(u)|$ is found
  for different values of $\sigma$, $x_0$, and confidence bound $r$. Top
  left: $x_0 = 1$, $r = 0.1$ and $\sigma \in \{0.01,0.02,\dots,0.17\}$;
  bottom left: $x_0 = 1$, $\sigma = 0.04$ and $r \in \{0.03, 0.06,\dots ,
  0.45\}$; and right: $x_0 = 1$ and $(\sigma, r) \in
  \{(0.01,0.03),\dots,(0.12,0.3)\}$. }
\end{figure}

\begin{figure}
  \centering
  \includegraphics[width=.24\textwidth,keepaspectratio]{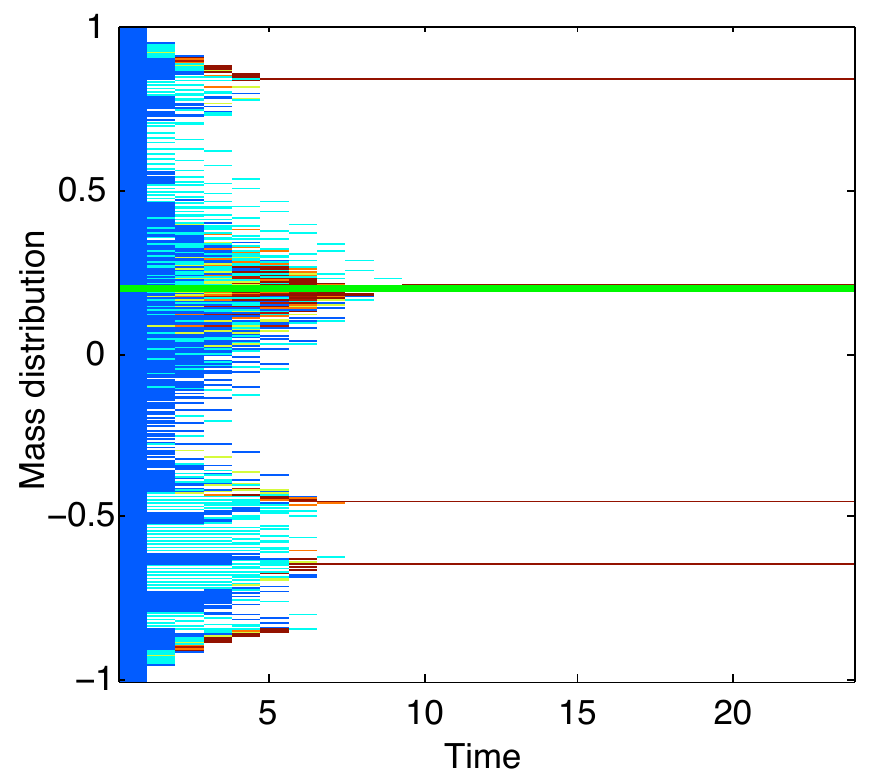}
  \includegraphics[width=.24\textwidth,keepaspectratio]{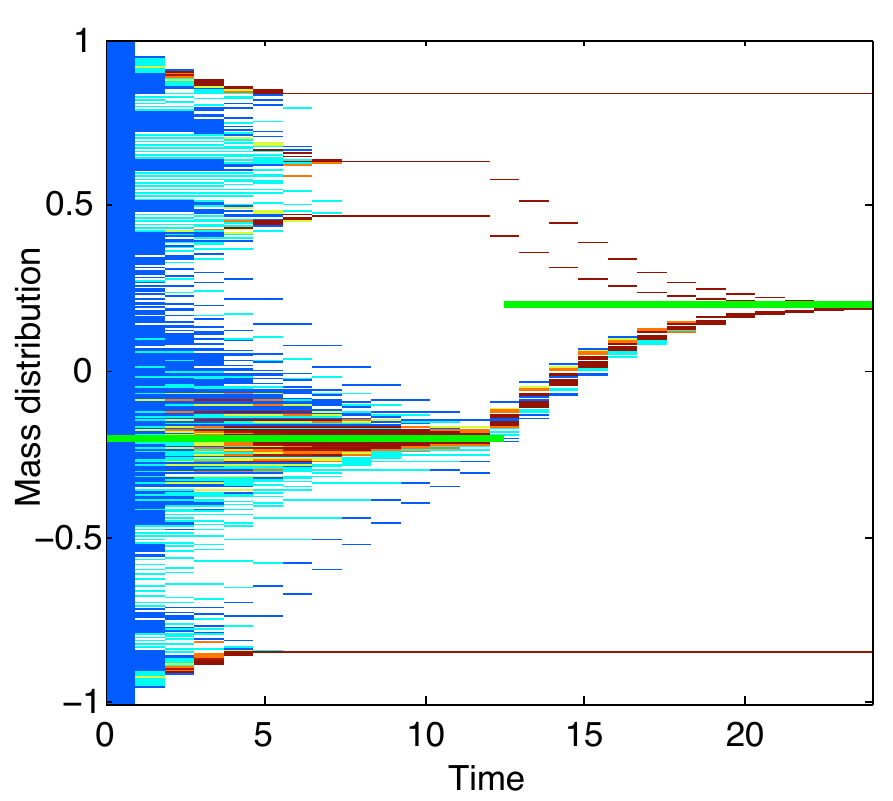}
  \caption{Two Eulerian HK systems with $\mu_0 \sim \mathcal{U}(-1,1)$ under the influence of direct strategy with $u \sim \subscr{\mathcal{N}}{truncated}(0.2,0.1^2)$ (left) and distracting strategy $u \sim \subscr{\mathcal{N}}{truncated}(-0.2,0.1^2)$ for $ t \le 12$ and $u \sim \subscr{\mathcal{N}}{truncated}(0.2,0.1^2)$ for $12 < t \le 25$ (right). In direct case 0.6525 portion of population is attracted to the input's center (green line), while in distracting case this portion is  0.8675. 
   }\label{figCompare}
\end{figure}

In this section we assume the exogenous input is a truncated Gaussian
distribution with bounded support centered at the advertised opinion, and
we  present some preliminary results on the influence of such an input
on the final population opinion value.
First, we study the effect of a single constant input.  We assume that
input is a truncated Gaussian distribution defined as follows: a Gaussian
distribution with mean $\hat{x}$ and variance $\sigma$ is truncated to have
support equal to $[\hat{x} - 3\sigma, \hat{x}+3\sigma]$. We denote such a
truncated Gaussian by $\subscr{\mathcal{N}}{truncated}(\hat{x},\sigma^2)$. The bounded
support assumption is based on the 3-sigma rule, where the effect of input
on our dynamical system outside the 3-sigma region can be ignored.
\begin{mydef}[Attraction range]
Consider an Eulerian HK system with input $u \sim
\subscr{\mathcal{N}}{truncated}(\hat{x},\sigma^2)$.  We define
\emph{attraction range} of input $u$, denoted by $\mathcal{R}(u)$, to be
the maximal opinion interval $[y,z] \in \real$ with the property that
$$ \lim_{t \rightarrow \infty} \gamma_t \circ \dots \circ \gamma_0 (y) =
\lim_{t \rightarrow \infty} \gamma_t \circ \dots \circ \gamma_0 (z) =
\hat{x}.$$
\end{mydef}
If the system satisfies conditions of Theorem~\ref{propertiesProp}, then
$\mu_0(\mathcal{R}(u))$ represents the \emph{attracted population}, i.e.,
the total population that reaches an opinion consensus at the center of the
input. The simulations reported in Figure~\ref{figLinearrelation} reveal a
linear relation between $|\mathcal{R}(u)|$, $\sigma$, and $r$ for
$|\mathcal{R}(u)| < 0.6 |\supp \mu_0|$ in the evolutions of an Eulerian HK
system with uniform initial distribution $\mu_0 \sim \mathcal{U}(-x_0,x_0)$
and input $u \sim \subscr{\mathcal{N}}{truncated}(0,\sigma^2)$. The
simulations lead us to an interesting conjecture.
\begin{conj}\label{conjLinearrelation}
Consider an Eulerian HK system with uniform initial mass distribution
$\mu_0$ whose support is a closed interval and input $u \sim
\subscr{\mathcal{N}}{truncated}(\hat{x},\sigma^2)$, where $\hat{x} \in
\supp \mu_0$.  If $\sigma$ is sufficiently small, then $|\mathcal{R}(u)| =
a \sigma + b r +c,$ with $a,b \in\real_{>0}$ and $c \in \real$.
\end{conj}
Second, in an Eulerian HK systems with $\mu_0 \sim \mathcal{U}(-x_0,x_0)$
and truncated Gaussian input, we study simple policies aimed at maximizing the
population with positive opinions in finite time $T$, or equivalently,
aimed at maximizing the objective function $\int_0^1 d\mu_T(z)$.  We
discuss two manipulation strategies$-$\emph{direct} and \emph{distracting}.
In direct strategy, the manipulator broadcasts a positive opinion for all
times. On the contrary, in distracting strategy the manipulator first
broadcasts a neutral or mildly negative opinion to attract the attention of
people with strong negative opinions, and only later broadcasts the
positive opinions.  Loosely speaking, the distracting strategy implements a
well-known persuation method: in dealing with someone with different
beliefs, a manipulator would start with a moderate opinion to win the trust
of that person. 
More precisely, if we assume that the input's variance is a fixed parameter
and the input is influential on less than half of entire population at
$t=0$ so that $\sigma < |\supp \mu_0|/12$, then we can define input's mean
based on our strategy.
\begin{itemize}
\item {\it Direct strategy:} Manipulator advertises for positive opinions,
  and thus $u_t \sim \subscr{\mathcal{N}}{truncated}(x_I(t),\sigma^2)$ with
  $0 < x_I(t) < |\mathcal{R}(u)|/2$ for all $t \le T$.
\item {\it Distracting strategy:} First, for all $ t\le \alpha T$ for some
  $\alpha \in (0,1)$, the manipulator broadcasts negative opinions, thus
  $u_t \sim \subscr{\mathcal{N}}{truncated}(x_{\ii}(t),\sigma^2)$ with
  $x_{\ii}(t) < 0$. One can assume that $x_{\ii}(t) = \xmin(t) +
  |\mathcal{R}(u)|/2$. Then, for all $\alpha T \le t \le T$, the
  manipulator advertises for positive opinions, thus $u_t \sim
  \subscr{\mathcal{N}}{truncated}(x_I(t),\sigma^2)$ with $0 < x_I(t) <
  |\mathcal{R}(u)|/2$.
\end{itemize}
The following discussion is a heuristic explanation for why the distracting
strategy outperforms the direct strategy.  It follows from the assumption
$\sigma < |\supp \mu_0|/12$ and boundedness of $|\mathcal{R}(u)|$ that the
direct strategy prevents attraction of the population with opinions in the
interval $[\xmin(0), -\mathcal{R}(u)/2]$.  However, in distracting
strategy, this population is in the attraction range of the first input
before $\alpha T$, and hence there is a fluctuation of population centered
at $x_{\ii}(t)$ and closer to the input's mean after $\alpha T$. An example of
this comparison is depicted in Figure~\ref{figCompare}.

 \section{Conclusion}\label{SecCon}

To describe the formation of opinions in a large population, we focused on an Eulerian model and introduced a reasonable exogenous input. 
First, we proved some fundamental properties of this dynamical system and derived a simple sufficient condition for consensus. Second, we proved the convergence of the mass distribution with no input to an atomic measure. Third, we defined the \emph{attraction range} of a normally distributed input, and for a uniformly distributed initial population over opinions, we conjectured a linear relation between attraction range's length and systems parameters. Accordingly, we compared two different manipulation strategies that aim to increase the population who vote positively in a finite time.

\bibliography{../ref/alias,../ref/Main,../ref/New}

\end{document}